\documentclass[leqno, 12pt]{article}
\usepackage{amssymb}

\def \EE {\mathbb{E}}
\def \RR {\mathbb{R}}

\def \e  {\epsilon}

\def \bx {{\bf x}}

\def \X  {\mathcal{X}}

\def \P  {\mathcal{P}}
\def \N  {\mathcal{N}}
\def \M  {\mathcal{M}}

\def \oon {\frac{1}{n}}

\def \olp  {\overline{p}}
\def \olen  {\overline{\epsilon_n}}
\def \oleni  {\overline{\epsilon_{ni}}}
\def \olenj  {\overline{\epsilon_{nj}}}

\newtheorem{proposition}{Proposition}[section]

\title{Asymptotic equivalence of the jackknife and infinitesimal jackknife variance estimators for some smooth statistics}
\author{Alex D. Gottlieb \\ 5745 Beck Ave, N Hollywood, CA  91601  U.S.A.}

\date{
 Running title: {\it Asymptotic equivalence of jackknives}
}

\begin{document}
\maketitle

\begin{abstract}

The jackknife variance estimator and the the infinitesimal
jackknife variance estimator are shown to be asymptotically
equivalent if the functional of interest is a smooth function of
the mean or a trimmed L-statistic with H\"older continuous weight
function.

\end{abstract}

\newpage

\section{Introduction}

\setcounter{equation}{0}

This note concerns the asymptotic behavior of the jackknife
variance estimator $v_{jack}$, especially regarding its
relationship to the infinitesimal jackknife variance estimator
$v_{ijack}$. We consider, in particular, the variance estimates
$v_{jack}$ and $v_{ijack}$ for smoothly trimmed L-statistics and
for smooth functions of the sample mean.  We prove that $v_{jack}$
and $v_{ijack}$ are often asymptotically equivalent to one another
in the sense that
\begin{equation}
\label{AsymptEquiv}
      v_{jack} -  v_{ijack}   \ = \    \mathrm{O}_{\mathrm{p}}\big(n^{-h}\big)
\end{equation}
for some $h > 0$.  The equivalence of $v_{jack}$ and $v_{ijack}$
can sometimes be used to prove that $v_{jack}$ is asymptotically
normal, but --- remarkably --- it also holds even when these
estimators are not asymptotically normal.

These issues are relevant to the following scenario of statistical
practice: One wishes to estimate some functional $T(p)$ of an
unknown population distribution $p$, and, to this end, one draws
$n$ samples from the population and uses $T(\epsilon_n)$ to
estimate $T(p)$, where $\epsilon_n$ is the empirical distribution
of the $n$ samples. One would then like an estimate of the
sampling variance of $T(\epsilon_n)$.   Two widely-used
nonparametric estimates of this variance are $v_{jack}$ and
$v_{boot}$, the jackknife and bootstrap variance estimates. Having
obtained one of these estimates, one naturally desires to know how
accurate it is.   To assess the accuracy of the usual Monte Carlo
approximation of $v_{boot}$, one may use the
jackknife-after-bootstrap technique of Efron (1992).   This note
concerns the asymptotic behavior of $v_{jack}$ and a closely
related estimator $v_{ijack}$, the infinitesimal jackknife.

Beran(1984) showed that $v_{jack}$, $v_{ijack}$, and $v_{boot}$
are asymptotically equivalent, and asymptotically normal, if the
functional $T(p)$ has a well-behaved second-order functional
derivative.  While the proof of Beran (1984) requires a strong
statement of the Dvoretsky-Kiefer-Wolfowitz inequality to handle
the asymptotics of $v_{boot}$, the asymptotic equivalence of
$v_{jack}$ and $v_{ijack}$ is easier to prove, and does not
require second-order differentiability of $T$.    Indeed, we shall
see that $v_{jack}$ and $v_{ijack}$ are asymptotically equivalent
if $T$ has a well-behaved first-order derivative and $v_{ijack}$
is consistent as an estimator of the variance of $T(\epsilon_n)$.

The equivalence of $v_{jack}$ to $v_{ijack}$ can help one to
determine the asymptotic variance of the former. For example,
Gardiner and Sen (1979) have carefully studied the asymptotic
normality of $v_{ijack}$ in the context of L-statistics.  We shall
see in Section~\ref{L} that their work establishes the asymptotic
normality of $v_{jack}$, too, thanks to the equivalence of
$v_{jack}$ and $v_{ijack}$ for variance estimation of
L-statistics.

Sometimes $v_{jack}$ and $v_{ijack}$ are asymptotically equivalent
even when they are not asymptotically normal.   Consider, for
example, the estimation of the variance of a function $g$ of the
sample mean.   If $g$ is once, but not twice, continuously
differentiable, then $v_{jack}$ and $v_{ijack}$ may not be
asymptotically normal, yet they will satisfy (\ref{AsymptEquiv})
as long as $g'$ is H\"older continuous of order $h$.

After the necessary definitions are presented in the next section,
we prove in Section~\ref{meanie} that $v_{jack}$ and $v_{ijack}$
are asymptotically equivalent as estimators of the variance of
smooth functions of the sample mean.  In Section~\ref{L} we
discuss the asymptotic normality of $v_{jack}$ as an estimator of
the variance of trimmed L-statistics.

\section{Background and definitions}

\setcounter{equation}{0}

Let $p$ be a probability measure on a sample space $\X$. Given $n$
samples from $\X$, sampled independently under the probability law
$p$, one desires to estimate the value $T(p)$ of some real
functional $T$ on the space $\P(\X)$ of all probability measures
on $\X$.  Denote by $\e_n$ the map that converts $n$ data points
$x_1,x_2,\ldots,x_n$ into the empirical measure
\begin{equation}
\label{e}
       \e_n(x_1,x_2,\ldots,x_n) \ = \ \oon \sum_{i=1}^n \delta(x_i)
\end{equation}
where $\delta(x_i)$ denotes a point-mass at $x_i$.  The {\it
plug-in estimate} of $T(p)$ given the data $\bx =
(x_1,\ldots,x_n)$ is
\begin{equation}
\label{plug-in}
         T_n \ = \  T(\e_n(\bx)).
\end{equation}
Suppose $T_n$ is an asymptotically normal estimator of $T(p)$,
i.e., suppose the distribution of $n^{1/2}(T_n - T(p))$ tends to
$\N(0,\sigma^2)$. The jackknife is a computational technique for
estimating $\sigma^2$: one transforms the $n$ original data points
into $n$ {\it pseudovalues} and computes the sample variance of
those pseudovalues.

Given the data $\bx = x_1,x_2,\ldots,x_n$, the {\it jackknife
pseudovalues} are
\[
         Q_{ni} \ = \ n T_n(\e_n) \ - \
(n-1)T(\e_{ni}) \qquad \qquad  i = 1,2,\ldots,n
\]
with $\e_n$ as in (\ref{e}) and
\begin{equation}
\label{e-ni}
     \e_{ni} \ = \ \frac{1}{n-1} \sum_{j\ne i} \delta(x_j).
\end{equation}
 The {\it jackknife variance
estimator} is
\begin{equation}
\label{JackVar}
           v_{jack}(x_1,x_2,\ldots,x_n) \ = \ \frac{1}{n-1}\sum_{i=1}^n
                     \left( Q_{ni} - \overline{Q_n} \ \right)^2
\end{equation}
where $\overline{Q_n} = \oon \sum Q_{nj}$.  The variance estimator
$v_{jack}$ is said to be {\it consistent} if
$v_{jack}\longrightarrow \sigma^2$ almost surely as $n \rightarrow
\infty$.  Sufficient conditions for the consistency of $v_{jack}$
are given in terms of the functional differentiability of $T$. An
early result of this kind states that $v_{jack}$ is consistent if
$T$ is strongly Fr\'echet differentiable (Parr(1985)), and it is now
known that $v_{jack}$ is consistent even if $T$ is only
continuously G\^ateaux differentiable as defined in Shao(1993).

A functional derivative of $T$ at $p$, denoted $\partial T_p$, is
a linear functional that best approximates the behavior of $T$
near $p$ in some sense.  For instance, a functional $T$ on the
space of bounded signed measures $\M(\X)$ is {\it G\^ateaux
differentiable} at $p$ if there exists a continuous linear
functional $\partial T_p$ on $\M(X)$ such that
\[
    \lim_{t \rightarrow 0} \big| t^{-1}\left(T(p + t m) - T(p)\right) \ - \   \partial T_p(m)
    \big| \ = \ 0
\]
for all $m \in \M(\X)$.   The concept of Hadamard
differentiability is more relevant to statistical asymptotics, for
the fluctuations of $T(\e_n)$ about $T(p)$ are asymptotically
normal if $T$ is Hadamard differentiable at $p$.  A functional
$T:\P(\RR) \longrightarrow \RR$ is {\it Hadamard differentiable}
at $p$ if there exists a continuous linear functional $\partial
T_p$ on $\M(\RR)$ such that
\[
    \lim_{t\rightarrow 0} \big|t^{-1} \left(T(p + t m_t) - T(p)\right)  \ - \   \partial T_p(m)
    \big| \ = \ 0
\]
whenever $\{m_t\}_{t \in \RR}$ is such that
$\lim\limits_{t\rightarrow 0} m_t = m$ and $m_t(\RR)=0$ for all
$t$, the topology on $\M(\RR)$ being the one induced by the norm
$\|m\| = \sup\limits_{t \in
\RR}\left\{\big|m((-\infty,t])\big|\right\}$.    If $T$ is
Hadamard differentiable at $p$,
 the variance of $n^{1/2}T(\e_n)$ tends to
\begin{equation}
\label{sigmasquared}
           \sigma^2 \ = \ \EE_p\phi_p^2
\end{equation}
as $n \longrightarrow \infty$, where $\phi_p(x)$ is the {\it
influence function}
\begin{equation}
\label{influence}
          \phi_p(x) \ = \  \partial T_p (\delta(x)- p)
\end{equation}
(this can be shown via the Delta Method using Donsker's theorem
(van der Waart(1998))). The {\it infinitesimal jackknife
estimator} (Jaeckel(1972)) of $\sigma^2$ is obtained by
substituting the empirical measure $\e_n$ for $p$ in
(\ref{sigmasquared}):
\begin{equation}
\label{infinitesimal}
         v_{ijack} \ = \   \EE_{\e_n} \phi_{\e_n}^2 .
\end{equation}

\section{Functions of the mean}
\label{meanie}

\setcounter{equation}{0}

When $q$ is a measure, we denote $\int x q(dx)$ by $\overline{q}$
if the integral is defined.  Let $g \in C^1(\RR)$ and let
\[
   T(m) \ = \ g \left(\overline{m}\right)
\]
be defined for all finite signed measures $m$ with finite first
moment. The functional derivative at $m$ of $T$, evaluated at $q$,
is $
     \partial T_m(q) = g'\left(\overline{m}
     \right)\overline{q}
$; the influence function (\ref{influence}) is $\phi_m(x) =
g'\left(\overline{m}
     \right) \left(x - \overline{m}\right)$.
Suppose that $x_1,x_2,\ldots$ are iid $p$, and $p$ has a finite
second moment. Let $T_n$ denote the plug-in estimator defined in
(\ref{plug-in}). Then the asymptotic variance of
$n^{1/2}\left(T_n-T(p)\right)$ is
\begin{equation}
\label{sigsquared}
         \sigma^2 \ = \ g'(\olp)^2\Big\{ \int x^2 p(dx) - \olp^2 \Big\}.
\end{equation}
Let $v_{jack}$ and $v_{ijack}$ denote the jackknife and
infinitesimal jackknife variance estimates of $\sigma^2$.

\begin{proposition}
\label{Main}
 If $g'$ is H\"older continuous of order
$h > 1/2$ (with global H\"older constant) and $p$ has a finite
moment of order $2+2h$ then
\[
      v_{jack} -  v_{ijack}
                \ = \
                  \mathrm{O}_{\mathrm{p}}\big(n^{-h}\big).
\]
\end{proposition}

\noindent {\bf Proof}: \qquad   Setting $\Delta_{ni} =
\left(Q_{nj} - \overline{Q_n}\right) - \phi_{\e_n}(x_i)$, formula
(\ref{JackVar}) for $v_{jack}$ yields
\begin{equation}
\label{GottaLabel}
    v_{jack}   \ = \
\EE_{\e_n}\phi_{\e_n}^2 \ + \ \frac{1}{n-1}\EE_{\e_n}\phi_{\e_n}^2
\ + \  \frac{2}{n-1} \sum_{i=1}^n \phi_{\e_n}(x_i)\Delta_{ni} \ +
\
                            \frac{1}{n-1} \sum_{i=1}^n
                            \Delta_{ni}^2 \ .
\end{equation}
The second term on the right hand side of (\ref{GottaLabel}) is $
\mathrm{O}_{\mathrm{p}}\left( 1/n \right)$ since
\[
      \EE_{\e_n}\phi_{\e_n}^2 \ = \
        \oon \sum_{i=1}^n \phi_{\e_n}^2(x_i) \ = \ \oon \sum_{i=1}^n
g'\left(\olen \right)^2(x_i - \olen)^2
\]
converges almost surely to $\sigma^2$.

 To control the last two terms on the right hand side of (\ref{GottaLabel})
 we need a bound on $\Delta_{ni}$.  Recall the
notation $\e_{ni}$ of (\ref{e-ni}).   Since $g$ is differentiable,
$ g \left( \olenj \right) - g \left(\oleni \right) \ = \
g'\left(\eta_{ji}\right)\left(\olenj - \oleni \right) $ for some
$\eta_{ji}$ between $\oleni$ and $\olenj$, so that
\[
     Q_{ni} - \overline{Q_n} \ = \
          \frac{n-1}{n} \sum_{j=1}^n
          \big(  g \left( \olenj \right) - g \left(\oleni \right)\big)
          \ = \
          \frac{n-1}{n} \sum_{j=1}^n
          g'\left(\eta_{ji}\right)\left(\olenj - \oleni \right).
\]
Therefore, since $\phi_{\e_n}(x_i) = g'\left(\olen \right)(x_i -
\olen) = \oon\sum_j g'\left(\olen \right)(x_i - x_j)$,
\begin{eqnarray*}
     \Delta_{ni}
     \ = \
     \left(Q_{ni} - \overline{Q_n}\right) - \phi_{\e_n}(x_i)
     & = &
        \frac{n-1}{n} \sum_{j=1}^n
          g'\left(\eta_{ji}\right)\left(\olenj - \oleni \right)
          \  - \  \oon\sum_{j=1}^n g'\left(\olen \right)(x_i - x_j) \\
     &  = &
                   \oon \sum_{j=1}^n
         \left( g'\left(\eta_{ji}\right) - g'\left(\olen \right)\right)(x_i - x_j
         ).
\end{eqnarray*}
But $g'$ is H\"older continuous of order $h$ and
$|\eta_{ji}-\olen|< \max\{|\olenj - \olen|,|\oleni - \olen|\}$, so
the identity $(n-1)\left( \e_n - \e_{ni}\right) \ = \
\delta_{x_i} - \e_n $ implies that
\[
            \left|  g'\left(\eta_{ji}\right) - g'\left(\olen
            \right)\right| \ \le \  C \big( |\olenj - \olen|^h +
            |\oleni - \olen|^h \big)
            \ \le \  C (n-1)^{-h} \big( |\olen - x_j|^h +
            |\olen - x_i|^h \big),
\]
where $C$ is a global H\"older constant for $g'$.  It follows that
\[
     \left| \Delta_{ni} \right|
     \ = \
     C (n-1)^{-h} \oon \sum_{j=1}^n   \big( |\olen - x_j|^h +
            |\olen - x_i|^h  \big) \big( |\olen - x_j| +
            |\olen - x_i|  \big).
\]
With this bound on $\Delta_{ni}$, and  assuming that $p$ has a
finite moment of order $2(1+h)$, it may be shown that
\[
     \oon \sum_{i=1}^n \Delta_{ni}^2 \ = \  \mathrm{O}_{\mathrm{p}}\big(n^{-2h}),
\]
and then, by the Cauchy-Schwartz inequality, that
\[
         \Big\vert \oon \sum_{i=1}^n \phi_{\e_n}(x_i)\Delta_{ni}
         \Big\vert \ = \  \mathrm{O}_{\mathrm{p}} \big( n^{-h} \big).
\]
Substituting the preceding estimates in (\ref{GottaLabel})
completes the proof. \hfill $\square$

Consider $g(x) = x - \hbox{sgn}(x)x^2$.  This function has a
Lipschitz continuous derivative but does not have a second-order
derivative at $0$.  By Proposition~\ref{Main},
\[
      v_{jack} -  v_{ijack}
                \ = \
                 \mathrm{O}_{\mathrm{p}}\left(1/n\right).
\]
However, for some population distributions having mean $0$, one
can prove that $v_{ijack}$ is not asymptotically normal, and
simulations suggest that $v_{jack} -  v_{boot}$ is
$\mathrm{O}_{\mathrm{p}}\left(1/\sqrt{n} \right)$ rather than
$\mathrm{O}_{\mathrm{p}}\left(1/n\right)$.  This example shows
that $v_{jack}$ is more closely related to $v_{ijack}$ than it is
to $v_{boot}$.

\section{Trimmed L-statistics} \label{L}

\setcounter{equation}{0}

Suppose that $\ell:(0,1) \longrightarrow \RR$ is supported on
$[\alpha,1-\alpha]$ for some $0 < \alpha < 1/2$, and let
\begin{equation}
\label{LFunctional}
     L(p) \ = \ \int_0^1 P^{-1}(s)\ell(s)ds.
\end{equation}
Here $P^{-1}$ denotes the quantile function for $p$, i.e.,
$
     P^{-1}(s) = \min\{x: P(x) \ge s\}
$ for $0<s<1$ where $P$ denotes the cdf of $p$.  A plug-in
estimate for $L$ is called a {\it trimmed L-statistic}, or a
trimmed {\it linear combination of quantiles}.

 If the weight function $\ell$ is
continuous then $L$ is Hadamard differentiable at all $p \in
\P(\RR)$ (see, e.g., Lemma 22.10 of van der Waart (1998)), and so
the L-statistics are asymptotically normal (an original reference
is Stigler (1974)).  The asymptotic variance $\sigma^2$ of the
L-statistics may be estimated by $v_{jack}$, which  converges
almost surely to $\sigma^2$ if $\ell$ is continuous (Parr(1985),
Shao and Tu (1995)).    The jackknife and infinitesimal jackknife
would seem to be the only nonparametric methods of consistent
variance estimation for L-statistics, aside from the bootstrap
(Shucany and Parr(1982)).

We turn now to the question of the asymptotic normality of
$v_{jack}$.  In this regard, a variant of the L-functional
(\ref{LFunctional}) has been treated in the literature, namely
\begin{equation}
\label{LFunctionalVariant}
    \mathcal{L}(p)= \int x \ell(P(x))p(dx)\ .
\end{equation}
If $P$ is continuous and strictly increasing then $\mathcal{L}$ of
(\ref{LFunctionalVariant}) equals $L$ of (\ref{LFunctional}).
Beran (1984) proves that $v_{jack}$ for $\mathcal{L}$ is
asymptotically normal
--- and so is $v_{boot}$ --- if $\ell$ is continuously
differentiable and $p$ has bounded support.  Section 2.2.3 of Shao
and Tu (1995) incorrectly claims that $v_{jack}$ for $\mathcal{L}$
is asymptotically normal if $\ell$ is H\"older continuous of order
greater than $1/2$, and it also wrongly claims that the asymptotic
variance equals $\hbox{Var}\big(\phi_p^2\big)$, where $\phi_p$ is
the influence function of $\mathcal{L}$.  A detailed discussion of
those errors is given in an unpublished technical report (Gottlieb
(2001)).
 Nevertheless, a reworking of Definition~2.6 and Theorem~2.7 in
Shao and Tu (1995) leads us to the following general proposition,
which will presently be applied to the case where the
$T(\epsilon_n)$ are L-statistics:
\begin{proposition}
\label{Abstract}
   Let $\epsilon_n$ denote the empirical distribution
of $n$ iid samples from $p$, and let $v_{jack}$ and $v_{ijack}$
denote the jackknife and infinitesimal jackknife estimates of the
variance of $T(\epsilon_n)$.

   Let $\| q' - q \|$ denote the supremum of the absolute value of the
difference between the cdf's of $q'$ and $q$.  Suppose that there
exist positive constants $C$, $\delta$, and $h$ such that
\begin{equation}
\label{expansion}
        T(q')  \ = \   T(q) + \partial_q T(q'-q) + \mathcal{R}(q',q)
\end{equation}
 for all $q',q$ with $\| q' - p \|,\| q - p \| < \delta$, where the
remainder $| \mathcal{R}(q',q) | \le C \| q' - q \|^{1+h} $. Then
\[
          v_{jack} - v_{ijack} \ = \
           \mathrm{O}_{\mathrm{p}}\big(n^{-h}\big)
\]
if $v_{ijack}$ is bounded in probability.
\end{proposition}
The straightforward proof of this proposition proceeds like the
proof of Proposition~\ref{Main} above, except that
(\ref{expansion}) is used to bound $\Delta_{ni}$ in
(\ref{GottaLabel}).

Now, let $\mathcal{L}$ be a trimmed L-functional of the form
(\ref{LFunctionalVariant}) whose weight function $\ell$ is
H\"older continuous of order $h$.  Upon integrating the right hand
side of (\ref{LFunctionalVariant}) by parts, it becomes easy to
verify that $\mathcal{L}$ admits the expansion (\ref{expansion})
near any $p$. Since $v_{ijack}$ converges almost surely,
Proposition~\ref{Abstract} implies that
  $v_{jack}$ and $v_{ijack}$ are asymptotically equivalent.

This equivalence allows us to conclude that $v_{jack}$ is
asymptotically normal for many L-functionals, for Gardiner and
Sen (1979) have found hypotheses that guarantee the asymptotic
normality of $v_{ijack}$ for generalized L-functionals of the form
(\ref{LFunctional}).  They begin by assuming that the cdf $P$ of
the population distribution is continuous. In this case
\[
      v_{ijack} \ = \    \EE_{\e_n}\phi_{\e_n}^2 \ = \ \int\int
         \ell(P_n(y)) \left[  P_n(y \wedge z) - P_n(y)P_n(z)
         \right]
\ell(P_n(z))dydz \ .
\]
In order to make contact with the work of Gardiner and Sen
(1979), let us suppose that $P$ is continuous and strictly
increasing. Their hypotheses are general enough to apply to
non-trimmed L-statistics, but too complicated to be repeated here.
Suffice it to say that their theorem applies under our current
assumptions that $\ell$ is trimmed and that $P$ has no jumps or
flats, if it is assumed in addition that $P$ does not have very
heavy tails and that $\ell$ is piecewise continuously
differentiable with H\"older continuity of order greater than
$1/2$ at the cusps. (Imagine, for example, a piecewise-linear
$\ell$ whose graph is shaped like a desert mesa; this is one of
the weight functions recommended in Stigler (1973) for smoothly
trimmed means.)  In these cases $v_{jack}$ is asymptotically
normal as well, by Proposition~\ref{Abstract}.

\section{Acknowledgments}

The author would like to thank Steve Evans for his advice and
encouragement, and also Rudolf Beran.  Support from the ESI and
WPI in Vienna and the START project {\it Nonlinear Schr\"odinger
and quantum Boltzmann equations} (FWF Y-137) is acknowledged.

\end{document}